\documentclass{amsart}

\usepackage{amsfonts}
\usepackage{amsmath}
\usepackage{amssymb}
\usepackage{amscd}
\usepackage{mathrsfs}  
\usepackage{amsbsy}
\usepackage{amsthm}
\usepackage{graphicx}
\usepackage{fancyhdr}
\usepackage{color}

\usepackage[OT2,OT1]{fontenc}  
\newcommand\cyr{%
\renewcommand\rmdefault{wncyr}%
\renewcommand\sfdefault{wncyss}%
\renewcommand\encodingdefault{OT2}%
\normalfont \selectfont} \DeclareTextFontCommand{\textcyr}{\cyr}
\input epsf 
\newcommand{\be}{\begin{equation}}
\newcommand{\ee}{\end{equation}}

\def\sideremark#1{\ifvmode\leavevmode\fi\vadjust{\vbox to0pt{\vss
 \hbox to 0pt{\hskip\hsize\hskip1em
\vbox{\hsize2cm\tiny\raggedright\pretolerance10000
 \noindent #1\hfill}\hss}\vbox to8pt{\vfil}\vss}}}%

\begin{document}

\title{ Operator-valued Frames on C*-Modules}
\author{Victor Kaftal, David Larson, Shuang Zhang* }
\address{Department of Mathematical Sciences, University of Cincinnati,
Cincinnati, Ohio 45221-0025}
\email{victor.kaftal@uc.edu }
\address{Department of Mathematics, Texas A{\&}M University, College Station, Texas 77843, USA}
\email{larson@math.tamu.edu}
\address{Department of Mathematical Sciences, University of Cincinnati,
Cincinnati, Ohio 45221-0025}
 \email{zhangs@math.uc.edu}


\maketitle
\date{}

\thanks{ * The first and third authors are partially supported by Taft Foundations, the second author is partially supported by an NSF grant. }
\begin{abstract}
Frames on Hilbert C*-modules have been defined for unital C*-algebras by Frank and Larson   [5] and operator-valued frames on a Hilbert space have been studied in [8]. The goal of this paper is to introduce operator-valued frames on a Hilbert C*-module for a $\sigma$-unital C*-algebra. Theorem 1.4 reformulates the definition given  in [5]  in terms of a series of rank-one operators converging in the strict topology. Theorem 2.2. shows that  the frame transform and the frame projection of an operator-valued frame are limits in the strict topology of a series in the multiplier algebra and hence belong to it. Theorem 3.3 shows that two
operator-valued frames are right similar if and only if they share the same frame projection. Theorem 3.4 establishes an one-to-one correspondence between Murray-von Neumann equivalence classes of projections in the multiplier algebra and right similarity equivalence classes of operator-valued frames and provides a parametrization of all Parseval operator-valued frames on a given Hilbert C*-module. Left similarity is then defined and Proposition 3.9 establishes when two left unitarily equivalent frames are also right unitarily equivalent.
\end{abstract}


\section* {Introduction}

Frames on a Hilbert space are collections of vectors satisfying the condition
$$a\Vert \xi\Vert ^2\le \sum_{j\in \mathbb J} \vert <\xi, \xi_j>\vert ^2 \le b\Vert \xi\Vert ^2$$ for some
positive constants $a$ and $b$ and all vectors $\xi$. This notion has been naturally extended by Frank and Larson  [5]  to countable collections of vectors  in a Hilbert C*-module for a unital C*-algebra satisfying an analogous defining property (see below {\bf 1.1} for the definitions). Most properties of frames on a Hilbert space hold also for Hilbert C*-modules, often have quite different proofs,  but new phenomena do arise.

A different generalization where frames are no longer vectors in a Hilbert space but operators on a Hilbert space is given in [8]  with the purpose of providing a  natural framework for multiframes, especially for those obtained from a unitary system, e.g, a discrete group representation. Operator-valued frames both generalize vector frames and can be decomposed into vector frames.

The goal of this article is to introduce the notion of operator-valued frame on a Hilbert C*-module.  Since the frame transform of Frank and Larson permits to identify a vector frame on an arbitrary Hilbert C*-module with a vector frame on the standard Hilbert C*-module $\ell^2(\mathcal A)$ of the associated C*-algebra $\mathcal A$, for simplicity's sake we confine our definition directly to frames on $\ell^2(\mathcal A).$  When the associated C*-algebra is $\sigma$-unital, it well known (see  [9] ) that the algebra of bounded adjointable operators on $\ell^2(\mathcal A)$ can be identified with the multiplier algebra $M(\mathcal A\otimes \mathcal K)$ of  $\mathcal A\otimes \mathcal K$,  about which a good deal is known. Since reference are mainly formulated in terms of right-modules, we treat  $\ell^2(\mathcal A)$ as a right module (i.e., as  `row vectors').

A frame on $\ell^2(\mathcal A)$ is thus defined as a collection of operators $\{A_j\}_{j\in \mathbb J}$ with \linebreak$A_j \in E_0M(\mathcal A\otimes \mathcal K)$ for a fixed projection $E_0\in M(\mathcal A\otimes \mathcal K)$ for which
$$aI\le \sum_{j\in \mathbb J} A_j^*A_j\le bI $$
for some positive constants $a$ and $b$, where $I$ is the identity of  $M(\mathcal A\otimes \mathcal K)$  and the series converges in the strict topology of $M(\mathcal A\otimes \mathcal K)$.

We will show in Theorem 1.4 how to associate (albeit not uniquely) to a vector frame in the sense of [5] an operator-valued frame. When $\mathcal A$ is unital, we will decompose in Section 3.10) every operator-valued frame (albeit not uniquely into vector frames (i.e., a multiframe)
Some properties of operator-valued frames on a Hilbert C*-module track fairly well the properties of operator-valued and vector-valued frames on a Hilbert space. Often, the key difference in the proofs is the need to express objects like the $frame$ $transform$ or the $frame$ $projection$ as  series of elements of $M(\mathcal A\otimes \mathcal K)$ that converge in the strict topology, and hence, belong to $M(\mathcal A\otimes \mathcal K)$.

We illustrate some commonalities  and differences with the Hilbert space case by considering in particular three topics. That the dilation approach of Han and Larson in  [6], which was extended to operator-valued frames on Hilbert spaces in [8], has a natural  analog for operator-valued frames on Hilbert C*-modules if the frame transform is defined to have values inside the same Hilbert C*-module instead of into an ampliation of it.

Similarity of frames can also be defined and characterized as in the Hilbert space case, but now there is also a similarity from the left and we compare the two notions.

Finally, there is a natural composition of operator-valued frames  -  a new operation that has no vector frame analog and that illustrates the `multiplicity' of operator-valued frames.

In this paper we have explored the analogs of some of the properties of Hilbert space frames and much work remains to be done.
\

\

\

\section{ Operator-valued frames }
\

\noindent {\bf 1.1. Frames and operator-valued frames on a Hilbert space.}

A frame on a Hilbert space $\mathcal H$ is a collection of vectors
$\{\xi_j\}_{j\in \mathbb J}$ indexed by a countable set $\mathbb J$ for which there exist
 two positive constants $a$ and $b$ such that for all $\xi\in \mathcal H$,Ê
$$a\Vert \xi\Vert ^2\le \sum_{j\in \mathbb J} \vert <\xi, \xi_j>\vert ^2 \le b\Vert \xi\Vert ^2.$$
Equivalently,
$$aI\le \sum_{j\in \mathbb J} \xi_j \otimes \xi_j\le bI, $$
where $I$ is the identity of $\mathcal B(\mathcal H)$, $\eta \otimes \xi$ is the rank-one operator defined by \linebreak
$(\eta \otimes \xi)\zeta:= <\zeta, \xi>\eta$ and
the series converges in the $strong$ $operator$ $topology$ (pointwise convergence). The above condition
 can be rewritten asÊ
$$aI\le \sum_{j\in \mathbb J} A_j^*A_j\le bI $$
where $A_j:= \eta \otimes \xi_j$ for some arbitrary fixed unit vector $\eta\in \mathcal H$. It is thus equivalent to  the series $\sum_{j\in \mathbb J} A_j^*A_j$ converging in the strong operator topology to a bounded invertible operator. Notice that the convergence of the numerical and the operatorial series are unconditional.

\

This reformulation naturally leads to the more general  notion of  operator-valued frames $\{A_j\}_{j\in \mathbb J}$  on a Hilbert space $\mathcal H$ in [8], namely a collection of operators \linebreak $A_j \in B(H, H_0)$, with ranges in a fixed Hilbert space $H_0$ (not necessarily of dimension one)  for which the series $\sum_{j\in \mathbb J} A_j^*A_j$ converges in the strong operator topology to a bounded invertible operator.  Operator-valued frames on a Hilbert space can be decomposed into, and hence identified with, multiframes.

\

Frames with values in a Hilbert C*-module have been introduced  in  [5] and then studied in  [7], [13],  and others. So much of the Hilbert space frame theory carries over, that one could argue that frame theory finds a natural general setting in Hilbert C*-modules. We will show in Theorem 1.4 that frames on a Hilbert C*-module can be equivalently defined in terms  of  rank-one operators on the module. This leads naturally to the definition of general operator-valued frames on a Hilbert C*-module. Before giving the formal definitions, we recall for the readers' convenience some relevant background about Hilbert C*-modules.
\

\

\noindent {\bf 1.2. Hilbert C*-modules} ([2,Ch.13], [9]).

\

  Let $\mathcal A$ be a  C*-algebra. Then a  Hilbert (right) C*-$\mathcal A$-module is a pair \linebreak
  $(\mathcal H, <.,.>)$,  with $\mathcal H$ a (right) module over  $\mathcal A$ and $<.,.>$ a binary operation from  $\mathcal H $ into $\mathcal A$, that
satisfies the following six axioms, similar to those of Hilbert spaces, except
that for right modules the linearity occurs for the second and not the first component of the  inner product.
 For $\xi,\eta , \eta _1, \eta _2\in \mathcal H$ and $a\in \mathcal A$
$$\aligned (i)& \ \ <\xi , \eta _1+\eta _2>=<\xi , \eta _1>+<\xi , \eta _2>;\\
(ii)&\ \ <\xi  , \eta a >=<\xi , \eta >a;\\
(iii)&\ \ <\xi ,\eta >^*=<\eta ,\xi>;\\
(iv)&\ \ <\xi ,\xi >\ge 0;\\
(v)&\ \ <\xi,\xi >=0\ \Longleftrightarrow \xi =0;\\
(vi)&\ \ (\mathcal H,\Vert \ .\ \Vert)   \text{ is complete, where} \Vert \xi \Vert := \Vert<\xi,\xi>\Vert^{1/2} \\
\endaligned $$
 
The  classic example of  Hilbert (right) $\mathcal A$-module and the only one we will consider in this paper is the $standard$ $module$ $\mathcal H_{\mathcal A}:=\ell^2(\mathcal A)$, the space of all sequences $\{a_i\}\subset \mathcal A$ such that
$\sum_{i=1}^\infty a_i^*a_i$ converges in norm to  a  positive element of $\mathcal A$. $\ell^2(\mathcal A)$ is endowed with the natural linear (A-module) structure and right $\mathcal A$-multiplication, and with
the $\mathcal A$-valued inner product  defined by
$$<\{a_i\},\{b_i\}>=\sum_{i=1}^\infty a_i^*b_i,$$
where the sum converges in norm by the Schwartz Inequality ([2] or [9]).

A  map $T$ from $\mathcal H_{\mathcal A}$ to $\mathcal H_{\mathcal A}$ is called a (linear) bounded operator on $\mathcal H_{\mathcal A}$, if  \linebreak $T(\lambda \xi  )+T(\mu \eta  ) = T(\lambda \xi  +\mu \eta  )$,     $T(\xi a ) = T(\xi)a$ for all $\xi, \eta \in \mathcal H_{\mathcal A}$, $ \lambda, \mu \in \mathbb C$, and $a\in \mathcal A$, and  if
 $$\Vert T\Vert :=\sup\{ \Vert T\xi \Vert \ \mid \xi \in \mathcal H_{\mathcal A}, \Vert \xi \Vert \le 1\} < \infty. $$
\noindent Not every bounded operator $T$ has a bounded adjoint $T^*$, namelyÊÊ
$$<T^* \xi ,\eta >=<\xi ,T\eta >\ \ \text{for all}\ \ \ \xi , \ \etaÊ
\  \in \mathcal H_{\mathcal A},$$
as there is no Riesz Representation Theorem for general Hilbert C*-modules.
Nevertheless, there are abundant operators on $\mathcal H_{\mathcal A}$ whose adjoints exist and the collection of bounded adjointable operators is denoted by  $ \mathcal B(\mathcal H_{\mathcal A})$.  Then,    $  \mathcal B(\mathcal H_{\mathcal A})$ is a C*-algebra (see [9] or [2, Ch. 13]).
Notice that if $\mathcal A=\mathbb C$, then  $\mathcal H_{\mathcal A}=\ell^2$ and $  \mathcal B(\mathcal H_{\mathcal A})=   \mathcal B(\ell^2)$. Some of the properties of   $ \mathcal B(\ell^2)$ extend naturally to
 $ \mathcal B(\mathcal H_{\mathcal A})$.
 For each pair of elements $\xi$ and $\eta $ in $\mathcal H_{\mathcal A}$,
 a bounded `rank-one'  operator is defined by
$$\theta _{\xi ,\eta }(\zeta)=\xi <\eta ,\zeta>\ \ \ \
\text{for all }  \zeta \in  \mathcal H_{\mathcal A}.$$
The closed linear span of all rank-one operators is denoted by $\mathcal K(\mathcal H_{\mathcal A})$. When  $\mathcal A=\mathbb C$,  $\mathcal K(\mathcal H_{\mathcal A})$ coincides with the ideal  $\mathcal K$ of  all compact operators on $\ell^2$.  $\mathcal K(\mathcal H_{\mathcal A})$ is always a closed ideal of $\mathcal B(\mathcal H_{\mathcal A})$, but contrary to the separable  infinite dimensional Hilbert space case, in general it is not unique (e.g., see [2] or [9]).

The  analog of the strong$^*$-topology on $\mathcal B(\ell^2)$ is the $strict$
$topology $  on  $\mathcal B(\mathcal H_{\mathcal A})$ defined by
$$\mathcal B(\mathcal H_{\mathcal A})\ni T_{\lambda }   \longrightarrow  T \ \text{strictly}\    \text{if } \Vert (T_{\lambda }-T)S\Vert \to 0 \text{ and  }
  \Vert S(T_{\lambda }-T)\Vert \to 0  \ \forall  S \in
\mathcal K(\mathcal H_{\mathcal A}).$$
 We will use the following elementary properties:   $T_{\lambda } \longrightarrow  T$ strictly iff \linebreak $T^*_{\lambda } \longrightarrow  T^*$ strictly, and either of these convergences implies $BT_{\lambda } \longrightarrow  BT$   and  $T_{\lambda }B \longrightarrow  TB$ strictly  for all
$B\in \mathcal B(\mathcal H_{\mathcal A})$. Also, if $T_{\lambda } \longrightarrow  T$ strictly and $S_{\lambda }  \longrightarrow  S$ strictly, then  $T_{\lambda }S_{\lambda } \longrightarrow TS $   strictly.

\

There is an alternative view of the objects $\mathcal B(\mathcal H_{\mathcal A})$ and $\mathcal K(\mathcal H_{\mathcal A})$. Embed the tensor product $\mathcal A\otimes \mathcal K$
into its Banach space double dual $(\mathcal A\otimes \mathcal K)^{**}$, which, as is well known,  is a W*-algebra ([12]). The multiplier algebra of $\mathcal A\otimes \mathcal K$, denoted by
$M(\mathcal A\otimes \mathcal K)$, is defined as the collection
$$\{T\in (\mathcal A\otimes \mathcal K)^{**}:\ \ TS, ST\in \mathcal A\otimes \mathcal K
\ \ \forall \ S\ \in \ \mathcal A\otimes \mathcal K\}.$$
 Equipped with the norm of
$(\mathcal A\otimes \mathcal K)^{**}$,ÊÊ
$M(\mathcal A\otimes \mathcal K)$ is a C*-algebra.  Assuming that
$\mathcal A$ is $\sigma$-unital, we will frequently apply the following two *-isomorphismsÊ
without further reference:Ê
$$ \mathcal B(\mathcal H_{\mathcal A})\cong M(\mathcal A\otimes \mathcal K)\ \ \ \ \text{and}\ \
\ \ \mathcal K(\mathcal H_{\mathcal A})\cong \mathcal A\otimes \mathcal K\
\   \ \ [\text{Ka1}]. $$ The algebra $\mathcal B(\mathcal H_{\mathcal A})$  is technicallyÊ
hard to work with, while $M(\mathcal A\otimes \mathcal K)$ is more accessible due to
many established results.  More information on the subject can be found in  the sample references
[9] and [2], among many others. Although most properties hold with appropriate modifications also for left modules, since the original theory was developed by Kasparov for right  Hilbert C*-modules ([9]), the results found in the literature are often formulated for right modules. This is the reason why  our definition of frames is given for right modules instead of left modules as in [5].Ê

To avoid unnecessary complications, from now on, we assume that $\mathcal A$ is a $\sigma$-unital C*-algebra.
\

\

\noindent {\bf 1.3. Vector Frames on Hilbert C*-modules}

According to  [5], a (vector) frame on the Hilbert C*-module $\mathcal H_{\mathcal A}$ of a  $\sigma$-unital C*-algebra  $\mathcal A$ is a collection of elements $\{\xi_j\}_{j\in \mathbb J}$ in $\mathcal H_{\mathcal A}$   for which there are two positive scalars $a$ and $b$ such that for all
$\xi\in\mathcal H_{\mathcal A}$,Ê
$$a<\xi,\xi>\le \sum_{j\in \mathbb J}  <\xi,\xi_j,><\xi_j,\xi> \le b<\xi,\xi>,$$
where the convergence is in the norm of the C*-algebra $\mathcal A$.Ê
The following theorem permits us to reformulate this definition in terms of rank-one operators.
Notice that $<\xi,\xi_j,><\xi_j,\xi>= <\theta_{\xi_j, \xi_j} \xi, \xi >$.

\

 {\bf 1.4. Theorem } Let $\mathcal A$ be a $\sigma$-unital  C*-algebra. Then the collection $\{\xi_j\}_{j\in \mathbb J}$ in the Hilbert C*-module $\mathcal H_{\mathcal A}$ is a frame if and only if
the series $\sum_{j\in \mathbb J}  \theta_{\xi_j,\xi_j}$ converges in the strict topology to a bounded  invertible operator in $B(\mathcal H_\mathcal A)$.

\

First we need the following elementary facts. For the readers' convenience we present their proofs.

\

 {\bf  1.5  Lemma } Assume that $\eta, \eta', \xi, \xi' \in \mathcal H_{\mathcal A}$. Then the following hold:
\begin{itemize}
 \item[(i)] $\theta _{\xi, \eta }\theta _{\eta ', \xi '}  =  \theta _{\xi <\eta, \eta'> ,\xi'}.$Ê
\item[(ii)] $\theta _{\xi ,\eta } ^* = \theta _{\eta  ,\xi }$.
\item[(iii)] $\theta _{\xi ,\eta }^*\theta _{\xi ,\eta }= \theta _{\eta < \xi ,\xi >, \eta }Ê
=\theta _{\eta < \xi ,\xi >^{\frac{1}2}, \eta <\xi , \xi > ^{1/2}}.$
\item[(iv)] If $T\in B(\mathcal H_\mathcal A)$, then $T\theta _{\xi ,\eta } = \theta _{T\xi ,\eta }$
\item[(v)] $\Vert \theta_{\xi, \eta}\Vert =\Vert \xi <\eta, \eta>^{1/2} \Vert  = \Vert <\xi, \xi>^{1/2}<\eta,\eta >^{1/2} \Vert$.  In particular, if  $\mathcal A$ is unital and $<\eta, \eta>= I$, then $\Vert \theta_{\xi, \eta}\Vert =\Vert \xi \Vert.$
\item[(vi)]$\Vert \theta_{\xi, \eta} \Vert  \le  \Vert \eta \Vert  \Vert \xi \Vert $
\item [(vii)] The rank-one operator $\theta _{\eta, \eta}$ is a projection if and only if $<\eta, \eta>$ is a projection, if and only if $\eta = \eta<\eta,\eta>$.
 \cite [Lemma 2.3]{5}
\end{itemize}
\begin{proof}
(i) For any $\gamma \in \mathcal H_{\mathcal A}$, one has
$$\aligned \theta _{\xi, \eta }\theta _{\eta ', \xi '} \gamma
&=\xi < \eta , \theta _{\eta ', \xi '} \gamma >\\
&= \xi  < \eta ,  \eta '< \xi ' ,  \gamma  >>\\
&=\xi  < \eta ,  \eta '>< \xi ' ,  \gamma  >\\
&=\theta _{\xi < \eta ,  \eta '>, \xi '} \gamma . \endaligned $$
(ii) $$\aligned <\theta _{\xi ,\eta } ^*\xi', \eta'>
&= <\xi', \theta _{\xi ,\eta } \eta'>\\
&=<\xi', \xi<\eta, \eta'>>\\
&=<\xi', \xi><\eta, \eta'>\\
&= <\eta<\xi, \xi'>, \eta'>\\
&= <\theta _{\eta  ,\xi }\xi', \eta'>. \endaligned $$
 
(iii) The first identity follows from (i) and (ii). Moreover
$$\aligned \theta _{\eta < \xi ,\xi >, \eta }\xi'
&=  \eta < \xi ,\xi ><\eta , \xi'  >\\
&=\eta < \xi ,\xi >^{\frac 12} (<\eta < \xi ,\xi >^{\frac 12}>, \xi'  > \\
&=  \theta _{\eta < \xi ,\xi >^{\frac{1}2}, \eta <\xi , \xi > ^{1/2}}\xi' . \endaligned $$

(iv) follows directly from the definition.

(v) $$\aligned \Vert \theta_{\xi, \eta}\Vert &=
\sup \{\Vert \theta_{\xi, \eta}\gamma\Vert \mid \Vert\gamma \Vert =1\}\\
&= \sup \{\Vert <\xi <\eta, \gamma>,\xi <\eta, \gamma>>\Vert^{1/2} \mid \Vert\gamma \Vert =1\}\\
&= \sup \{\Vert <\eta, \gamma>^*<\xi ,\xi><\eta, \gamma>\Vert^{1/2} \mid \Vert\gamma \Vert =1\}\\
&= \sup \{\Vert <\xi ,\xi>^{1/2}<\eta, \gamma>\Vert\mid \Vert\gamma \Vert =1\}\\
&= \sup \{\Vert <\eta<\xi ,\xi>^{1/2}, \gamma>\Vert\mid \Vert\gamma \Vert =1\}\\
&= \Vert \eta<\xi ,\xi>^{1/2}\Vert\\
&=\Vert <\xi ,\xi>^{1/2}<\eta,\eta><\xi ,\xi>^{1/2}\Vert^{1/2}\\
&=\Vert <\xi ,\xi>^{1/2}<\eta,\eta>^{1/2}\Vert
\endaligned $$

(vi) is obvious

(vii) For completeness we add a short proof.  By (ii), $\theta _{\eta, \eta}$  is a projection if and only if
$$0 = \theta _{\eta, \eta}-\theta _{\eta, \eta}\theta _{\eta, \eta}= \theta _{\eta, \eta} - \theta _{\eta<\eta,\eta>, \eta} = \theta _{\eta-\eta<\eta,\eta>, \eta}$$Ê
and by (v), this condition is equivalent to $(\eta-\eta<\eta,\eta>)<\eta,\eta>^{1/2}=0$.
If  $<\eta,\eta>$ is a projection, then
 $$<\eta-\eta<\eta, \eta>,\eta-\eta<\eta, \eta>>= <\eta, \eta>- 2<\eta, \eta>^2 + <\eta, \eta>^3=0,$$
hence $<\eta-\eta<\eta, \eta>=0$, and thus $\theta _{\eta, \eta}$  is a projection. Conversely,
if \linebreak $(\eta-\eta<\eta,\eta>)<\eta,\eta>^{1/2}=0$ then
$$<\eta, (\eta-\eta<\eta,\eta>)<\eta,\eta>^{1/2}> = <\eta,\eta>^{3/2}-<\eta,\eta>^{5/2}=0,$$
whence $<\eta, \eta>$ is a projection.

\end{proof}

 \
 
 Notice that  equality in (vi) may fail.  For instance, if $\xi = \{p, 0, 0, ..., \}$ and \linebreak
 $\eta = \{q, 0, 0, ..., \}$ where $p,~q\in \mathcal A$ are orthogonal non-zero projections, then \linebreak
 $<\xi,\xi> = p= <\xi ,\xi>^{1/2}$ and $<\eta,\eta>=q=<\eta,\eta>^{1/2}$ hence  \linebreak $<\xi ,\xi>^{1/2}<\eta,\eta>^{1/2}=0$.

\

\begin{proof} [Proof of Theorem 1.4]
Assume first that  the series $\sum_{j\in \mathbb J}  \theta_{\xi_j,\xi_j}$ converges in the strict topology to some operator $D_A\in B( \mathcal H_{\mathcal A})$. Set $T_{\mathbb F} =\sum_{j\in \mathbb F}\theta_{\xi_j,\xi_j} -D_A$ for any finite subset $\mathbb F$ of $\mathbb J$.  Then the  net
$\{ T_{\mathbb F} \}$  converges to 0 in the strict topology.
Using the equality $\Vert T_{\mathbb F} \theta_{\xi, \eta}\Vert = \Vert T_{\mathbb F}\xi <\eta, \eta >^{1/2} \Vert $ for all
$\xi, \eta$ from Lemma 1.5 (v), it follows that $\Vert T_{\mathbb F}\xi a \Vert \to 0$ for all positive $a\in \mathcal A\otimes \mathcal K$. But then,
  $$\aligned
&\Vert T_{\mathbb F}\xi \Vert   \le \Vert T_{\mathbb F}\xi a\Vert + \Vert T_{\mathbb F}(\xi -a\xi)\Vert \\
& \le \Vert T_{\mathbb F}\xi a\Vert + \sup\{\Vert T_{\mathbb F}\Vert \Vert(\xi -a\xi)\Vert\\
&= \Vert T_{\mathbb F}\xi a\Vert + \sup\{\Vert T_{\mathbb F}\Vert \Vert\{<\xi, \xi>-a<\xi, \xi> -<\xi, \xi>a +a<\xi, \xi>a)\}\Vert^{1/2}\\
&\le \Vert T_{\mathbb F}\xi a\Vert + \sup\{\Vert T_{\mathbb F}\Vert \{\Vert(<\xi, \xi> -a<\xi, \xi>)\Vert + \Vert a\Vert \Vert<\xi, \xi> -a<\xi, \xi>\Vert\}^{1/2}.
\endaligned$$
Since every C*-algebra $\mathcal A$ has a positive  approximate identity, one can  choose $a > 0$ such that  $$\sup\{\Vert T_{\mathbb F}\Vert \{\Vert(<\xi, \xi> -a<\xi, \xi>)\Vert + \Vert a\Vert \Vert<\xi, \xi> -a<\xi, \xi>\Vert\}^{1/2} < \epsilon.$$
For that $a$, $\Vert T_{\mathbb F}\xi a\Vert  < \epsilon$ for all $\mathbb F \supset G$ for some finite subset $G$ of $\mathbb J$. This shows that  $\Vert T_{\mathbb F}\xi \Vert \rightarrow 0 $.
Consequently,
 the series  $\sum_{j\in \mathbb J}  \theta_{\xi_j,\xi_j}\xi$ converges in the norm of $\mathcal H_{\mathcal A}$ to
$D_A\xi$, and hence,  $$\sum_{j\in \mathbb J}  <\theta_{\xi_j,\xi_j}\xi, \xi> =  \sum_{j\in \mathbb J}  <\xi,\xi_j,><\xi_j,\xi>$$ converges in the norm of $\mathcal A$ to $<D_A\xi, \xi> $ by the Schwartz Inequality. Now a positive operator $D_A$ is bounded and invertible if and only if $aI \le D_A \le bI$ \linebreak  for some constants $a,~b > 0$. By \cite[2.1.3] {11}, this condition is equivalent to  \linebreak
$a<\xi,\xi>~ \le <D_A\xi, \xi> ~\le b<\xi,\xi>.$
 Therefore, $\{\xi_j\}_{j\in \mathbb J}$ is a frame.

Conversely, assume that $\{\xi_j\}_{j\in \mathbb J}$ is a frame. Then
  $\sum_{j\in \mathbb J}  <\theta_{\xi_j,\xi_j}\xi, \xi>$ converges in the norm of $\mathcal A$ to $<D_A\xi, \xi>$ for some positive operator $D_A\in B( \mathcal H_{\mathcal A})$.
For any finite subset  $\mathbb F \subset \mathbb J$, $<\sum_{j\in \mathbb F}  \theta_{\xi_j,\xi_j}  \xi, \xi>
\le <D_A\xi, \xi>,$ hence, again by \cite[2.1.3]{11},
$$0 \le D_A - \sum_{j\in \mathbb F}  \theta_{\xi_j,\xi_j} \le D_A.$$
But then, by (vi) in the above lemma,
$$\aligned\Vert (D_A - \sum_{j\in \mathbb F}  \theta_{\xi_j,\xi_j})\theta_{\xi,\eta}\Vert &\le
\Vert (D_A - \sum_{j\in \mathbb F}  \theta_{\xi_j,\xi_j})^{1/2}\Vert  \Vert \theta_{(D_A - \sum_{j\in \mathbb F}  \theta_{\xi_j,\xi_j})^{1/2}\xi,\eta}\Vert\\
& \le \Vert D_A \Vert ^{1/2}  \Vert (D_A - \sum_{j\in \mathbb F}  \theta_{\xi_j,\xi_j})^{1/2}\xi \Vert \Vert \eta\Vert\\
&= \Vert D_A \Vert ^{1/2} \Vert \eta\Vert  \Vert <(D_A - \sum_{j\in \mathbb F}  \theta_{\xi_j,\xi_j})\xi , \xi>\Vert^{1/2}\to 0.
\endaligned$$
Since the linear span of rank-one operators is by dense in $\mathcal A \otimes \mathcal K$, it follows that $\Vert (D_A - \sum_{j\in \mathbb F}  \theta_{\xi_j,\xi_j})S\Vert \to 0$ for all $S\in \mathcal A \otimes \mathcal K.$  Since $D_A - \sum_{j\in \mathbb F}  \theta_{\xi_j,\xi_j}$ is selfadjoint, this proves that the series $\sum_{j\in \mathbb J}  \theta_{\xi_j,\xi_j}$ converges to $D_A$ in the strict topology.
The same argument as above shows that since  $D_A$ is bounded and invertible then
$$a<\xi,\xi>\ \le <D_A\xi, \xi>  =  \sum_{j\in \mathbb J}  <\xi,\xi_j,><\xi_j,\xi> \le b<\xi,\xi> \text{ for all } ~\xi.$$
\end{proof}
\
\

Many of the results on frames in Hilbert C*-modules are obtained under the assumption  that $\mathcal  A$ is unital, which is of course the case for Hilbert space frames where $\mathcal A=\mathbb C$. When   $\mathcal  A$ is unital,  in lieu of viewing frames as collections of vectors in $\mathcal H_{\mathcal A}$, we can view them as collections of rank-one operators on $\mathcal H_\mathcal A$ with range in a submodule $\mathcal H_0$. Indeed, if $\eta \in \mathcal H_{\mathcal A}$ is an arbitrary  unital vector, i.e., $<\eta, \eta>=I$, then by Lemma 1.5 (vii), (i) and (ii), $E_0:=\theta_{\eta, \eta}  \in \mathcal A\otimes K$ is a projection, actually the range projection of  $\theta_{\eta, \xi}$ for every $\xi  \in \mathcal H_{\mathcal A}$. Then  $\mathcal H_0:=E_0 \mathcal H_{\mathcal A}$ is a submodule of $\mathcal H_{\mathcal A}$ and we can identify $E_0M(\mathcal A\otimes \mathcal K)$ with $B(\mathcal H_\mathcal A,\mathcal H_0)$, the set of linear bounded adjointable operators from   $\mathcal H_\mathcal A$ to the submodule $\mathcal H_0$. Notice that all this would hold also under the weaker hypothesis that $<\eta, \eta>$ is a projection.
Then for every collection $\{\xi_j\}_{j\in \mathbb J}$ in  $\mathcal H_{\mathcal A}$,  define the rank-one operators $A_j:=\theta_{\eta, \xi_j}$. Since $E_0A_j=A_j$,  $A_j \in B(\mathcal H_\mathcal A,\mathcal H_0)$. Again, by Lemma 1.5, \linebreak$A_j^*A_j= \theta_{\xi_j<\eta, \eta>, \xi_j}=\theta_{\xi_j, \xi_j}$. It follows from Theorem 1.4  that $\{\xi_j\}_{j\in \mathbb J}$ is a frame if and only if  the series $ \sum_{j\in \mathbb J} A_j^*A_j $ converges  in the strict topology to a bounded invertible operator on $\mathcal H_\mathcal A$ .

This leads naturally to the following definition.

\

 \noindent {\bf 1.6 Definition  } Let $\mathcal A$   be a $\sigma$-unital C*-algebra  and $\mathbb J$ be a countable
index set.  Let $E_0$ be a projection in $M(\mathcal A\otimes \mathcal K)$. Denote by $\mathcal H_0 $ the submodule $ E_0\mathcal H_\mathcal A$  and identify  $B(\mathcal H_\mathcal A,\mathcal H_0)$ with $E_0M(\mathcal A\otimes \mathcal K)$.  A collection  $A_j\in B(\mathcal H_\mathcal A, \mathcal H_0)$ for $j\in \mathbb J$ is called an operator-valued frame on $\mathcal H_\mathcal A$ with range in  $\mathcal H_0$ if the sum $ \sum_{j\in \mathbb J} A_j^*A_j $ converges  in the strict topology to a bounded invertible operator on $\mathcal H_\mathcal A$, denoted by $D_A$.Ê
 $\{ A_j\}_{ j\in \mathbb J}$ is called a tight operator-valued frame (resp.,  a Parseval operator-valued frame) if $D_A= \lambda I $ for a positive number $\lambda$ (resp.,  $D_A=  I $).  If the set $\bigcup \{ A_j \mathcal H_\mathcal A: j\in \mathbb J\}$ is dense in $\mathcal H_0$, then the frame is said to be non-degenerate.

 From now on,  by frame, we will mean an operator-valued frame on a Hilbert C*-module. Notice that if  $\{ A_j\}_{ j\in \mathbb J}$ is a frame with range in $\mathcal H_0$ then it is also a frame with range in any larger submodule.
 
 A minor difference with the definition given in \cite{8} for operator-valued frames on a Hilbert space, is that  here, in order to avoid introducing maps between different modules,  we take  directly $\mathcal H_0$ as a submodule of $\mathcal H_\mathcal A$. For operator-valued frames on a Hilbert space we do not assume in [8]  that $\mathcal H_0\subset \mathcal H$ and hence we are left with the flexibility of considering $\dim \mathcal H_0  >  \dim \mathcal H$.

 \

\
 
  \noindent {\bf 1.7 Example } Let $\sum_{j\in \mathbb J}  E_j = I$ be a decomposition of the identity of $M(\mathcal A\otimes \mathcal K)$ \linebreak  into mutually orthogonal equivalent projections in $M(\mathcal A\otimes \mathcal K)$, i.e.,  $L_jL_j^*=E_j$ and $L_j^*L_j=E_0$ for some collection of partial isometries $L_j\in M(\mathcal A\otimes \mathcal K)$, and the series converges in the strict topology.  Let  $T$ be a left-invertible element of $M(\mathcal A\otimes \mathcal K)$ and let  $A_j:=L_j^*T$. Then $A_j \in B(\mathcal H_\mathcal A, E_0\mathcal H_\mathcal A) $ for $j\in \mathbb J$, and \linebreak$\sum_{j\in \mathbb J} A_j^*A_j=\sum_{j\in \mathbb J} T^*A_jT=T^*T$ is an invertible element of $M(\mathcal A\otimes \mathcal K)$, where the convergence is in the strict topology. Thus   $\{ A_j\}_{ j\in \mathbb J}$ is a frame with range in  $ E_0\mathcal H_\mathcal A$. The frame is Parseval precisely when  $T$ is an isometry. We will see in the next section that this example is generic.
\

\

\

\section{Frame Transforms }

\

\noindent {\bf 2.1. Definition }
Assume that $\{A_j\}_{ j\in \mathbb J}$ is a frame in $B(\mathcal H_\mathcal A, E_0\mathcal H_\mathcal A)$ for the Hilbert $C^*$-module $\mathcal H_\mathcal A$ and set $\mathcal H_0:=E_0\mathcal H_\mathcal A$. Decompose the identity ofÊÊ
$M(\mathcal A\otimes \mathcal K)$,
into a strictly converging sum of mutually orthogonal projections  $\{E_j\}_{j\in \mathbb J}$ in $M(\mathcal A\otimes \mathcal K)$ with
$E_j \sim E_{00} \ge E_0$. Let $L_j$ be partial isometries in $M(\mathcal A\otimes \mathcal K)$ such that  $L_jL_j^* = E_j$ and $L_j^*L_j = E_{00}$. Define the frame transform $\theta _A$ of the frame  $\{A_j\}_{ j\in \mathbb J}$ as
$$\theta _A=\sum_{j\in \mathbb J} L_j A_j:\ \ \mathcal H_\mathcal A\longrightarrow  \mathcal H_\mathcal A.$$
\

\

 {\bf  2.2. Theorem} Assume that  $\{A_j\}_{ j\in \mathbb J}$ is a frame in $B(\mathcal H_\mathcal A, \mathcal H_0)$.

\noindent {\bf (a) } The  sum $\sum_{j\in \mathbb J} L_j A_j$   converges in the strict topology, and hence $\theta _A$ is an element of $M(\mathcal A\otimes \mathcal K)$.

\noindent {\bf (b) }  $D_A= \theta_A^*\theta _A $,  $\theta_AD_A^{-1/2 }$ is an isometry, $P_A:=\theta_A D_A^{-1} \theta _A^*$ is the range projection of $\theta_A$, andÊ
all these three elements belong to $ M(\mathcal A\otimes \mathcal K).$

\noindent {\bf (c) }  $\{A_j\}_{ j\in \mathbb J}$ is a Parseval frame, if and only if $\theta _A$ is an isometry
of $M(\mathcal A\otimes \mathcal K)$,  and again  if and only if $\theta_A\theta _A^*$ is a projection.

\noindent {\bf (d) }  $A_j = L_J^*\theta_A$ for all $j\in \mathbb J.$

\begin{proof} {\bf (a) }  For every finite subset $F$ of $\mathbb J$, let $S_F=\sum_{j\in F} L_j A_j$.
We need to show that $\{S_F:F\ \text{is a finite subset of }\  \mathbb J\}$ is a Cauchy
net in the strict topology of $M(\mathcal A\otimes \mathcal K)$, i.e., for every for any $ a\in \mathcal A\otimes \mathcal K$, $\text{max}\{\Vert (S_F-S_{F'})a\Vert , \Vert a(S_F-S_{F'})\Vert  \longrightarrow 0$, in the sense that  for every $\varepsilon > 0$ there is a finite set $G$ such thatÊ
$$\text{max}\{\Vert (S_F-S_{F'})a\Vert , \Vert a(S_F-S_{F'})\Vert
< \varepsilon  \ \ \ \text{for any finite sets } F\supset G,~F' \supset G.$$Ê
Firstly, since the partial isometries $L_j$ have mutually orthogonal ranges, one has
$$   \aligned \Vert (S_F-S_{F'})a\Vert  &=
\Vert a^*(S_F-S_{F'})^*(S_F-S_{F'}) a\Vert ^{\frac{1}2}\\
&=\Vert a^*(\sum_{j\in (F  \setminus  F' )\cup (F'  \setminus  F)} A_j^*L_j^*L_jA_j) a\Vert ^{\frac{1}2}\\
&=\Vert a^*(\sum_{j\in (F  \setminus  F' )\cup (F'  \setminus  F)}A_j^*E_{00}A_j) a\Vert ^{\frac{1}2}\\
&\le \Vert a\Vert ^{\frac{1}2}\Vert \sum_{j\in (F  \setminus  F' )\cup (F'  \setminus  F)}A_j^*A_j
a\Vert ^{\frac{1}2}   \longrightarrow 0, \endaligned$$
where the last term above converges to 0 because of the assumption that $\sum_{j\in \mathbb J} A_j^*A_j $  converges in the  strict topology.
Secondly,  for all $\xi \in \mathcal H_\mathcal A$

$$ \aligned\Vert (S_F-S_{F'})\xi \Vert  ^2 &= \sum_{j\in (F  \setminus  F' )\cup (F'  \setminus  F)}\Vert L_jA_j\xi \Vert ^2\\
&= \sum_{j\in (F  \setminus  F' )\cup (F'  \setminus  F)}\Vert A_j\xi \Vert ^2\\
&= \sum_{j\in (F  \setminus  F' )\cup (F'  \setminus  F)}<A_j^*A_j\xi , \xi>\\
&=< \sum_{j\in (F  \setminus  F' )\cup (F'  \setminus  F)}A_j^*A_j\xi , \xi>\\
& \le <D_A \xi, \xi> \\
&= \Vert D_A^{1/2}\xi\Vert^2
\endaligned$$
Thus $\Vert S_F-S_{F'} \Vert \le \Vert D_A \Vert ^{1/2}$ for any finite sets $F$ and $F'$.
Moreover,  \linebreak
$ \sum_{j\in (F  \setminus  F' )\cup (F'  \setminus  F)}E_j(S_F-S_{F'}) = S_F-S_{F'}$, hence for every
$a\in \mathcal A\otimes \mathcal K$,

$$   \aligned \Vert a(S_F &-S_{F'}) \VertÊ
= \Vert a (S_F-S_{F'}) (S_F-S_{F'})^* a^*\Vert ^{\frac{1}2}\\
&= \Vert a  \sum_{j\in (F  \setminus  F' )\cup (F'  \setminus  F)}E_j(S_F-S_{F'}) (S_F-S_{F'})^*  \sum_{j\in (F  \setminus  F' )\cup (F'  \setminus  F)}E_ja^*\Vert ^{\frac{1}2}\\
&\le \Vert S_F-S_{F'} \Vert\Vert a  \sum_{j\in (F  \setminus  F' )\cup (F'  \setminus  F)}E_j a^*\Vert ^{\frac{1}2}\\
&\le \Vert D_A \Vert^{1/2} \Vert a \Vert ^{1/2}~\Vert \sum_{j\in (F  \setminus  F' )\cup (F'  \setminus  F)}E_j a^*\Vert^{1/2}\longrightarrow 0
\endaligned$$
by the strict convergence of the series $\sum_{j\in \mathbb J} E_j$.
The statements (b) and (c) are now obvious.

(d) Also obvious since  the series
$\sum_{i\in \mathbb J} L_j^*L_i A_i$ converges strictly to $L_j^*\theta_A$ and  $$L_j^*L_iA_i=\delta_{i,j}E_{00}A_i=\delta_{i,j}A_j.$$

\end{proof}
 
 The projection $P_A$ is called the {\sl frame projection}.
 \
 
 \
 
\noindent {\bf 2.3  Remark }
\noindent {\bf (i)} Since $\theta_A$ is left-invertible as  $(D_A^{-1}\theta_A^*)\theta_A = I$,  Example 1.7 is indeed generic, i.e., every frame $\{A_j\}_{j\in \mathbb J}$ is obtained from partial isometries $L_j$ with mutually orthogonal range projections  summing  to the identity and same first projection majorizing $E_0$.  The relation with the  by now familiar ``dilation" point of view of the  theory of frames is clarified in (ii) below.

\noindent {\bf (ii)} For vector  frames on a Hilbert space, the frame transform is generally defined as a map from the Hilbert space $\mathcal H$ into $\ell^2(\mathbb J)$ - a dilation of  $\mathcal H$.  If  $\mathcal H$ is infinite dimensional and separable and if $ \mathbb J$ is infinite and countable, which are the most common assumptions, then $\mathcal H$ can be identified with $\ell^2(\mathbb J)$, and hence, the frame transform can be seen as mapping of $\mathcal H$ onto a subspace. In the case of Hilbert C*-modules, it is  convenient to choose the latter approach, so to identify the frame transform and the frame projection with elements of $M(\mathcal A\otimes \mathcal K)$.

\noindent {\bf (iii)} The range projections of elements of $M(\mathcal A\otimes \mathcal K)$ and even of elements of $\mathcal A\otimes \mathcal K$ always belong  to $(\mathcal A\otimes \mathcal K)^{**}$, but may fail to belong to $M(\mathcal A\otimes \mathcal K)$. As shown above, however, the frame projection $P_A$ is always in $ M(\mathcal A\otimes \mathcal K)$  and $P_A \sim I$ since $\theta_A D_A^{1/2}$ is an isometry.

\noindent {\bf (iv)} When $\mathcal A$ is not simple, given an arbitrary nonzero projection \linebreak $E_0 \in M(\mathcal A\otimes \mathcal K)$ there may be no decomposition of the identity in projections equivalent to $E_0$.  Nevertheless,  there is always a decomposition of the identity into a strictly convergent sum of mutually orthogonal projections that are all equivalent to a projection $E_{00} \ge E_0$, e.g.,  $E_{00} = I$.

\noindent {\bf (v)} There seem to be no major advantage in considering only non-degenerate frames, i.e., seeking a ``minimal" Hilbert module $\mathcal H_0$ that contains the ranges of all  the operators $A_j$ or similarly,  choosing a frame transform with  a minimal  projection $E_{00}$. In fact, if we view the operators  $A_j$  as having their ranges in some $E_{00}\mathcal H_{\mathcal A}$, the ensuing frame transform will,  as seen in (d) above and in the next section, carry equally well all the ``information" of the frame.

\noindent {\bf (vi)} For the vector case, \cite [4.1] {5} proves that the frame transform $\theta$ as a map from a finite or countably generated Hilbert C*- module $\mathcal H$ to the standard module $\mathcal H_{\mathcal A}$ is adjointable.  This is obvious in the case that we consider, where $\mathcal H= \mathcal H_{\mathcal A}$ as  then $\theta_A$ is in the C*-algebra $M(\mathcal A\otimes \mathcal K)$.

\noindent {\bf (vii)} A compact form of the $reconstruction$ $formula$ for a frame is simply
$$D_A^{-1}\sum_{j\in \mathbb J} A_j^*A_j= D_A^{-1}\theta _A^*\theta _A = I.$$  In the special case that $\mathcal A$ is unital and that $\{A_j\}_{j\in \mathbb J}$ is a vector frame,  we have seen in the course of the proof of Theorem 1.4 that $\sum_{j\in \mathbb J} A_j^*A_j$ converges strongly and the same result was obtained in \cite[4.1]{5}. Thus, assuming for the sake of simplicity that the frame is Parseval,  the reconstruction formula has the more familiar form
$$\sum_{j\in \mathbb J} A_j^*A_j\xi = \sum_{j\in \mathbb J} \xi_j<\xi_j,\xi> = \xi \text{ for all } \xi \in \mathcal H_{\mathcal A}$$ where the convergence is in the norm of $\mathcal H_{\mathcal A}$.

\

\

\section{Similarity of frames }
\

\

\noindent {\bf  3.1. Definition } Two frames $\{A_j\}_{j\in \mathbb J}$  and $\{B_j\}_{j\in \mathbb J}$ in $B(\mathcal H_\mathcal A, \mathcal H_0)$ are said to be right-similar (resp. right-unitarily equivalent)
 if there exists an invertible (resp. a unitary) element
$T\in M(\mathcal A\otimes \mathcal K)$ such that $B_j=A_jT$ for all $j\in \mathbb J$.

\

The following facts are immediate and their proofs are left to the reader.

\

   {\bf 3.2.  Lemma }ÊÊ

\noindent {\bf (i)} If $\{A_j\}_{j\in \mathbb J}$ is a frame and $T$ is an invertible element in $M(\mathcal A\otimes \mathcal K)$, then $ \{A_jT\}_{ j\in \mathbb J}$ is also a frame.

 \noindent {\bf (ii) } If $\{A_j\}_{j\in \mathbb J}$  and $\{B_j\}_{j\in \mathbb J}$ are right-similar and
 $T$  is an invertible element in  $M(\mathcal A\otimes \mathcal K)$ for which $B_j=A_jT$ for all $j\in \mathbb J$, then $\theta _B=\theta _AT$. Therefore \linebreak $T= D_A^{-1}\theta _A^*\theta _B$, hence $T$ is uniquely determined. Moreover,  $P_A=P_B$ and \linebreak $D_B=T^*D_AT$. Conversely,
if $\theta _B=\theta_AT$ for some invertible element $T\in M(\mathcal A\otimes \mathcal K)$,
then $B_j=A_jT$ for all $j\in \mathbb J$.

\noindent {\bf (iii) } Every frame is right-similar to a Parseval frame, i.e., $\{A_j\}_{j\in \mathbb J}$  is right-similar to $\{A_j\}_{j\in \mathbb J}D_A^{-1/2}$. Two Parseval frame are right-similar if and only if they are right-unitarily equivalent.

\

 {\bf 3.3. Theorem  } Let $\{A_j\}_{j\in \mathbb J}$ and $\{B_j\}_{j\in \mathbb J}$ be two frames
in $B(\mathcal H_\mathcal A, \mathcal H_0)$. Then the following are equivalent:
\begin{itemize}
 \item[(i)] $\{A_j\}_{j\in \mathbb J}$ and $\{B_j\}_{j\in \mathbb J}$ are right-similar.
\item[(ii)] $\theta_B=\theta_AT$ for some invertible operator $T\in M(\mathcal A\otimes \mathcal K)$.
\item[(iii)] $P_A=P_B$.
\end{itemize}
\

\begin{proof} The implications  $(i)\Rightarrow (ii) \Rightarrow (iii)$ are given by Lemma 3.2. \linebreak
$(iii)\Rightarrow(i)$. Assume that $P_A=P_B$. Then by Theorem 2.2 (b) $$ \theta _A D_A^{-1}\theta _A^*=\theta _B D_B^{-1}\theta _B^*.$$
By Theorem 2.2. (d)Ê
$$B_j = L_j^*\theta _B = L_j^*P_A\theta _B =L_j^*\theta _AD_A^{-1}\theta _A^*\theta _B=
A_jD_A^{-1}\theta _A^*\theta _B.$$
Let $T:=D_A^{-1}\theta _A^*\theta _B$, then $T \in M(\mathcal A\otimes \mathcal K)$ and $$(D_B^{-\frac{1}2}\theta _B^*\theta _AD_A^{-\frac{1}2})(D_A^{-\frac{1}2}\theta _A^*\theta _BD_B^{-\frac{1}2})=
 D_B^{-\frac{1}2}\theta _B^*P_B\theta _BD_B^{-\frac{1}2} =I.$$
Interchanging $A$ and $B$, one has also
$$(D_A^{-\frac{1}2}\theta _A^*\theta _BD_B^{-\frac{1}2})(D_B^{-\frac{1}2}\theta _B^*\theta _AD_A^{-\frac{1}2})=I.$$
Thus, $D_A^{-\frac{1}2}\theta _A^*\theta _BD_B^{-\frac{1}2}$ is  unitary, hence $\theta _A^*\theta _B$ is invertible, and thus so is \linebreak $T=D_A^{-1}\theta _A^*\theta _B$. This concludes the proof.
 \end{proof}

\
As is the case in $B(H)$, all the projections in $M(\mathcal A\otimes \mathcal K)$ that are equivalent to the frame projection of a given frame, are also the frame projection of a frame, which by Theorem 3.3 is unique up to right similarity.\

\

 {\bf  3.4. Theorem  } Let $\{A_j\}_{j\in \mathbb J}$ be a frame in $B(\mathcal H_\mathcal A, \mathcal H_0)$ and let $P$ be a projection in $M(\mathcal A\otimes \mathcal K)$. Then $P\sim P_A$ if and only if there exists a frame $\{B_j\}_{j\in \mathbb J}$ in $B(\mathcal H_\mathcal A, \mathcal H_0)$ such that $P=P_B$.

\begin{proof} If $P=P_B$, let $V=\theta _BD_B^{-\frac{1}2} D_A^{-\frac{1}2} \theta _A$. Then $V\in B(\mathcal H_{\mathcal A})$, $VV^*=P$, and $V^*V=P_A$, i.e.,  $P\sim P_A$.
Conversely, if there exists $V\in M(\mathcal A\otimes \mathcal K)$ such that $VV^*=P$ and $V^*V=P_A$, then set $B_j=L_j^*V\theta_A$.
Then $$\sum_{j\in \mathbb J}B_j^*B_j = \theta_A^*V^*(\sum_{j\in \mathbb J}L_j^*L_j)V\theta _A=\theta_A^*V^* V\theta _A=D_A.$$
Thus  $\{B_j\}_{j\in \mathbb J}$ is a frame with $D_A=D_B$. Moreover, $$\theta _B =\sum_{j\in \mathbb J}L_j^*L_jV\theta _A = V\theta _A.$$
It follows that $P_B=VV^*=P.$
\end{proof}

The proof of Theorem 3.4 actually yields a parametrization of all the operator-valued frames with range  in $B(\mathcal H_\mathcal A, \mathcal H_0)$. For simplicity's sake, we formulate it in terms of Parseval frames
\

\

{\bf  3.5. Corollary} Let $\{A_j\}_{j\in \mathbb J}$ be a Parseval frame in $B(\mathcal H_\mathcal A, \mathcal H_0)$. Then $\{B_j\}_{j\in \mathbb J}$ is a Parseval frame in $B(\mathcal H_\mathcal A, \mathcal H_0)$ if and only if
$B_j=L_j^*V\theta_A$ for some partial isometry  $V\in M(\mathcal A\otimes \mathcal K)$ such that $VV^*=P_B$ and $V^*V=P_A$.

\

\

\noindent {\bf  3.6. Remark }
 For a given equivalence $P\sim P_A$  the choice of partial isometry
$V \in M(\mathcal A\otimes \mathcal K)$
with $VV^*=P $ and $V^*V=P_A$ is determined up to a unitary that commutes with $P_A$, i.e.,
 $V_1V_1^*=P$ and $V_1^*V_1=P_A$ implies $V_1=VU$ for some unitary $U$ that commutes with$P_A$, or, equivalently, $V_1=U_1V$ for some unitary $U_1$ that commutes with $P$.
Notice that if $V$ and $U$ are as above, then $\{L_j^*V\theta _AD_A^{-\frac{1}2}\}_{ j\in \mathbb J}$ and
$\{L_j^*VU\theta _AD_A^{-\frac{1}2}\}_{ j\in \mathbb J}$ are two frames having the same frame projections $P$.
It then follows from Theorem 3.3 that the two frames are right-similar, actually,
right-unitarily equivalent, because $L_j^*VU\theta _AD_A^{-\frac{1}2}=L_j^*V\theta _A
D_A^{-\frac{1}2}(D_A^{-\frac{1}2}\theta _A^*U\theta _AD_A^{-\frac{1}2})$
and $D_A^{-\frac{1}2}\theta _A^*U\theta _AD_A^{-\frac{1}2}$ is a unitary operator on $\mathcal H_\mathcal A$ as $\theta _AD_A^{-\frac{1}2}$ is an isometry.

\

\

For operator valued frames it is natural to consider also the notion of left similarity.

\

\noindent {\bf  3.7.  Definition } Two frames  $\{A_j\}_{j\in \mathbb J}$  and  $\{B_j\}_{j\in \mathbb J}$  in $B(\mathcal H_\mathcal A, \mathcal H_0)$
are said to be left-similar (resp., left-unitarily equivalent ) if there exists an invertible (resp., a unitary) element
$S$ in the corner algebra $E_0M(\mathcal A\otimes \mathcal K)E_0$ such that $B_j=SA_j$ for all $j\in \mathbb J$,
where $E_0$ is the projection in $M(\mathcal A\otimes \mathcal K) $ such that
$E_0\mathcal H_A=\mathcal H_0$.
\

\
The following results are elementary.

\

 {\bf  3.8. Lemma } If $\{A_j\}_{j\in \mathbb J}$  is a frame in $B(\mathcal H_\mathcal A, \mathcal H_0)$
and $S$ is an invertible element in $E_0M(\mathcal A\otimes \mathcal K)E_0$, then
$\{A_j\}_{j\in \mathbb J}$   is also a frame. Moreover, $\theta _B=\sum_{j\in \mathbb J} L_j SA_j$, hence $D_B=\sum_{j\in \mathbb J} A_j^*S^*SA_j$ and thus  $$\Vert S^{-1}\Vert ^{-2} D_A\le D_B\le \Vert S\Vert ^2D_A.$$
In particular, if $S$ is unitary, then $D_A=D_B$.

\

\

 {\bf 3.9. Proposition  } Given two left-unitarily equivalent frames $\{A_j\}_{j\in \mathbb J} $
and $\{B_j\}_{j\in \mathbb J} $  in $B(\mathcal H_\mathcal A, \mathcal H_0)$ with $B_j=SA_j$ for some unitary element
$S\in M(\mathcal A\otimes \mathcal K) $, then the following are equivalent:
\begin{itemize}
 \item[(i)] $S$ commutes with $A_jD_A^{-1}A_i^*$ for all $i,~j\in \mathbb J.$
\item[(ii)] $\{A_j\}_{j\in \mathbb J} $
and $\{B_j\}_{j\in \mathbb J} $ are right-unitarily equivalent.
\end{itemize}
 
\begin{proof}$(i) \Longrightarrow (ii)$ One has
$$\aligned  B_j=&SA_j\\
= &SA_j (D_A^{-1}D_A) =  \sum _{i\in \mathbb J} SA_jD_A^{-1}A_i^*A_i\\
=&\sum_{i\in \mathbb J} A_jD_A^{-1}A_i^*SA_i = A_jD_A^{-1} \sum_{i\in \mathbb J} A_i^*SA_i
\endaligned $$
where the series here and below converge in the strict topology.
Let $$U=D_A^{-\frac{1}2}\sum_{i\in \mathbb J} A_i^*SA_i  D_A^{-\frac{1}2}.$$ Then
$$\aligned UU^*=&D_A^{-\frac{1}2}\sum_{i,j\in \mathbb J} A_i^*SA_i  D_A^{-1} A_j^*S^*A_j  D_A^{-\frac{1}2} \\
=& D_A^{-\frac{1}2}\sum_{ i,j\in \mathbb J} A_i^*A_i  D_A^{-1} A_j^*SS^*A_j  D_A^{-\frac{1}2} \\
=& D_A^{-\frac{1}2}(\sum_{i\in \mathbb J} A_i^*A_i ) D_A^{-1}(\sum_{j\in \mathbb J} A_j^*A_j)  D_A^{-\frac{1}2}\\
=&I. \endaligned $$
Similarly,  $U^*U=I.$
Notice that   $B_j=A_jD_A^{-\frac{1}2}UD_A^{\frac{1}2}$, and that $D_A^{-\frac{1}2}UD_A^{\frac{1}2}$ is invertible.
It follows that $\{A_j\}_{j\in \mathbb J} $ and $\{B_j\}_{j\in \mathbb J} $ are right-unitarily equivalent frames.

$(ii) \Longrightarrow (i)$ $P_A=P_B$ by Theorem 3.3 and $D_A=D_B$ by Lemma 3.8. Then
$$\aligned P_A &=\theta _A D_A^{-1}\theta _A^* =  \sum_{i,j\in \mathbb J} L_iA_iD_A^{-1}A_j^* L_j^*\\
&=P_B=\theta _B D_B^{-1}\theta _B^* =  \sum_{i,j\in \mathbb J} L_iSA_iD_A^{-1}A_j^*S^* L_j^*.\endaligned$$
Multiplying on the left by $L_i^*$ and on the right by $L_j$, one has
$$A_iD_A^{-1}A_j^*=SA_iD_A^{-1}A_j^*S^*\ \ \ \forall\ i,j\in \mathbb J.$$Ê

\end{proof}

\noindent {\bf 3.10. Composition of frames } Let  $\{A_j\}_{\in \mathbb J} $ be a frame in $B(\mathcal H_\mathcal A,\mathcal H_0)$ and $\{B_i\}_{i\in \mathbb I}$ be a frame in $B(\mathcal H_0, \mathcal H_1)$. Then it is easy to check that $\{C_{i,j}:= B_iA_j\}_{j\in \mathbb J, i\in \mathbb I}$ is a frame in $B(\mathcal H_\mathcal A, \mathcal H_1)$,  called the composition of  the frames $\{A_j\}_{\in \mathbb J} $ and $B(\mathcal H_\mathcal A,\mathcal H_0)$. In symbols, $C=BA$
\

It is easy to see that the composition of two Parseval frames is also Parseval.\

\

\noindent {\bf 3.11. Remarks  } {\bf (i) } If $BA=BA'$, then $A=A'$. Indeed, if for all $i\in \mathbb I$ and $j\in \mathbb J$ and
$B_iA_j=B_iA'_j$, then $\sum_{i\in I} B_i^*B_iA_j=D_BA_j=D_BA'_j$. Then $A_j=A_j'$, since $D_B$ is invertible.
\

\

\noindent {\bf (ii) } If $A:= \{A_j\}_{\in \mathbb J} $ is non-degenerate (i.e., the closure of $\bigcup _{j\in \mathbb J} A_j\mathcal H_2$ is $\mathcal H_2$), then $BA=B'A$ implies $B=B'$.
\

\

\noindent {\bf 3.12. Remark } Let $\mathcal A$ be unital, then as shown in Theorem 1.4, we can identify a vector frame  $\{\xi_i\}_{i\in \mathbb I}$ on on a submodule $\mathcal H_1\subset \mathcal H_\mathcal A$, with the  (rank-one) operator valued frame $\{\theta_{\eta, \xi_i}\}_{i\in \mathbb I}$ in $B(\mathcal H_1, \theta_{\eta, \eta}\mathcal H_\mathcal A)$, where $\eta\in\mathcal H_\mathcal A$ is a vector for which $<\eta, \eta>= I$.  But then, for any operator-valued frame $\{A_j\}_{\in \mathbb J} $  in $B(\mathcal H_\mathcal A,\mathcal H_0)$ and $\mathcal H_1\subset \mathcal H_0$, the composition $\{\theta_{\eta, \xi_i}A_j= \theta_{\eta, A_j^*\xi_i}\}_{i \in \mathbb I,i \in \mathbb I}$ is identified with the vector frame $\{A_j^* \xi_i\}_{i \in \mathbb I,~j \in \mathbb J}$. We can view this as a decomposition of the operator valued frame $\{A_j\}_{\in \mathbb J} $ into the collection of (vector-valued) frames $\{A_j^* \xi_i\}_{j \in \mathbb J}$ indexed by $\mathbb I$, i.e., a ``multiframe".                                                                                        

\
\

\end{document}